\documentclass{llncs}

\usepackage{amscd,amsfonts,amsmath,amssymb,cite}
\usepackage{algorithm,algpseudocode}
\usepackage{graphicx}
\usepackage[colorlinks,linkcolor=blue,citecolor=blue,filecolor=blue,urlcolor=blue]{hyperref}

\begin{document}

\title{Computing Rational Generating Function of a Solution to the Initial Value Problem of Two-dimensional Difference Equation with Constant Coefficients}

\author{A.~A.~Kytmanov\inst{1} \and A.~P.~Lyapin\inst{2} \and T.~M.~Sadykov\inst{3}}

\institute{Siberian Federal University \\ Svobodny, 79,
Krasnoyarsk, 660041, Russia \\ \email{aakytm@gmail.com}
\and
Friedrich Schiller University of Jena \\ 07737 Jena, Germany \\
\email{LyapinAP@yandex.ru }
\and
Plekhanov Russian University \\
125993 Moscow, Russia \\ \email{Sadykov.TM@rea.ru}
}

\maketitle

\begin{abstract}
Algorithms for computing rational generating functions of
solutions of one-dimensional difference equations are well-known
and easy to implement. We propose an algorithm for computing
rational generating functions of solutions of two-dimensional
difference equations in terms of initial data of the corresponding
initial value problems. The crucial part of the algorithm is the
reconstruction of infinite one-dimensional initial data on the
basis of finite input data. The proposed technique can be used for
the development of similar algorithms in higher dimensions. We
furnish examples of implementation of the proposed algorithm.
\end{abstract}

\section{Introduction}

The study of linear difference equations and their solutions is an important topic of modern computer algebra (see, e.g., \cite{ABHP,AGH}).

The method of generating functions is a powerful tool for the study of differential equations with applications in the theory of discrete dynamical systems and in the enumerative combinatorial analysis. It allows one to apply the methods of complex analysis to the problems of enumerative combinatorics.

The one-dimensional case is well-studied (see \cite{gelfond,Isaacs}) due to the absence of geometric obstacles. In \cite{moivre1724}, A.~Moivre considered the power series
\[
f(0) + f(1) z + \ldots + f(k) z^k + \ldots
\]
with the coefficients $f(0), f(1),\ldots$ satisfying the difference equation
\begin{multline}\label{n1_1}
c_m f(x+m) + c_{m-1} f(x+m-1) + \ldots + c_0 f(x) = 0, \quad x=0,1,2\ldots,
\end{multline}
where $c_m \neq 0$, and $c_j \in \mathbb C$ are constants. He proved that this series always represents a rational function (De Moivre's Theorem, \cite{moivre1724}).

In the multidimensional case, which is not adequately explored (see \cite{petrovsek,lein2004,lein2007,LPZ2008}), the rational generating functions are the most useful class of generating functions according to the Stanley's hierarchy (see \cite{stanley1990}). A broad class of two-dimensional sequences that lead to rational generating functions is well-known in the enumerative combinatorics. For example, one can consider the
problem of finding a number of lattice paths, the problem on generating trees with marked labels, Bloom's strings, a number of placement of the pieces on the chessboard etc. (see \cite{Egorychev,merlini2008,bloom}).

Generating functions for multiple sequences with elements which could be expressed in terms of rational, exponential functions or gamma function, form a wide class of hypergeometric-type functions~\cite{Sadykov,DickensteinSadykovDoklady}. Their study leads to the problem of solving overdetermined systems of linear equations with polynomial coefficients.

A multidimensional analogue of the De Moivre Theorem was formulated and proved in \cite{LeinLyapin2009}. We now give some definitions and notations that we will need for formulating the main result.

\section{Known Results and Formulation of the Problem}

Let $x=(x_1,\ldots, x_n)\in \mathbb{Z}^n$, where $\mathbb{Z}^n = \mathbb Z \times \ldots \times \mathbb Z$ is the $n$-dimensional integer lattice. Let $A=\{\alpha\}$ be a finite set of points in $\mathbb{Z}^n$. Let $f(x)$ be a function of integer arguments $x=(x_1, \ldots, x_n)$ with constant coefficients $c_\alpha$. By a difference equation with respect to the unknown function $f(x)$ we call the equation of the form
\begin{equation}\label{razn_ur}
\sum_{\alpha \in A} c_\alpha f(x+\alpha) =0.
\end{equation}

In the present work we consider the case when the set $A$ belongs to the positive octant $\mathbb{Z}_0^n=\{(x_1,\ldots,x_n)\colon x_i\in\mathbb{Z}, x_i\geqslant 0, i=1,\ldots,n\}$ of the integer lattice and satisfies to the condition:
\begin{multline}\label{neravenstva}
\text{There exists a point}\ m=(m_1,\ldots, m_n)\in A \ \text{such that} \\ \text{for any}\ \alpha \in A \ \text{the inequalities}\ \alpha_j \leqslant m_j,\ j= 1,2,\ldots,n \ \text{hold}.
\end{multline}

Denote the \textit{characteristic polynomial} of (\ref{razn_ur}) by
\[
P(z) = \sum\limits_{\alpha\in A} c_\alpha z^\alpha = \sum\limits_{\alpha\in A} c_{\alpha_1,\ldots,\alpha_n} z_1^{\alpha_1} \cdots z_n^{\alpha_n}.
\]

The generating function (or the $z$-transformation) of the function $f(x)$ in integer variables $x\in \mathbb Z_0^n$ is defined as follows:
\[
F(z) = \sum_{x\geqslant 0} \frac {f(x)}{z^{x+I}}, \quad \text{where}\ I=(1,\ldots,1).
\]

For example, a solution of the difference equation
\[
f(x_1+1, x_2+1) - f(x_1+1,x_2) - f(x_1, x_2+1) + f(x_1, x_2)=0
\]
is an arbitrary function of the form $f(x_1,x_2) = \varphi(x_1)+\psi(x_2)$, where $\varphi$ and $\psi$ are arbitrary univariate functions of integer arguments and the corresponding generating series is in general divergent.

Define the ``initial data set'' for the difference equation (\ref{razn_ur}) that satisfies the condition (\ref{neravenstva}), as follows:
\[
X_0 = \{\tau \in \mathbb Z^n\colon \tau \geqslant 0,\ \tau \ngeqslant m\}.
\]
Here $\ngeqslant$ means that the point $\tau$ belongs to the complement of the set defined by the system of inequalities
\[
\tau_j \geqslant m_j,\ j=1,\ldots,n.
\]

The initial value problem is set up as follows: we need to find the solution $f(x)$ of the difference equation (\ref{razn_ur}), which coincides with a given function $\varphi(x)$ on $X_0$:
\begin{equation}\label{nach_dan}
f(x) = \varphi(x), \; x\in X_0.
\end{equation}

It is easy to show (see, e.g., \cite{lein2007}), that if the condition (\ref{neravenstva}) holds, then the initial value problem (\ref{razn_ur}), (\ref{nach_dan}) has the unique solution. The solvability of the problem (\ref{razn_ur}), (\ref{nach_dan}) without the constraints (\ref{neravenstva}) has been studied in \cite{petrovsek}.

We now proceed with some more notations that we will need later on.

Let $J=(j_1,\ldots,j_n)$, where $j_k \in \{0,1\}$, $k=1,\ldots, n$, be an ordered set. With each such set $J$ we associate the face $\Gamma_J$ of the $n$-dimensional integer parallelepiped
\[
\Pi_m = \{x\in\mathbb Z^n: 0\leqslant x_k \leqslant m_k, k=1,\ldots,n\}
\]
as follows:
\begin{align*}
\Gamma_J = \{x\in \Pi_m\colon &x_k = m_k,\ \text{if}\ j_k=1, \\
\text{and}\ &x_k < m_k,\ \text{if}\ j_k =0\}.
\end{align*}
For example, $\Gamma_{(1,\ldots,1)} = \{m\}$ and
\[
\Gamma_{(0,\ldots,0)}=\{x\in \mathbb Z^n\colon 0 \leqslant x_k < m_k,\ k=1,\ldots,n\}.
\]

It is easy to check that $\Pi_m = \bigcup\limits_J \Gamma_J$ and for any $J, J'$ the corresponding faces do not intersect: $\Gamma_J\cap \Gamma_{J'} = \emptyset$.

Let $\Phi(z) = \sum\limits_{\substack{\tau \geqslant 0 \\ \tau \ngeqslant m}} \dfrac{\varphi(\tau)}{z^{\tau+I}}$ be the generating function of the initial data of the solution of (\ref{razn_ur}), (\ref{nach_dan}). With each point $\tau\in\Gamma_J$ we associate the series
\[
\Phi_{\tau, J} (z) = \sum\limits_{y\geqslant \; 0} \frac{\varphi(\tau+Jy)}{z^{\tau+Jy+I}},
\]
and with each face $\Gamma_J$ we associate the series
\[
\Phi_J (z) = \sum\limits_{\tau \in \; \Gamma_J} \Phi_{\tau, J}(z).
\]
If we extend the domain of $\varphi(x)$ by zero on $\mathbb Z^n_+\setminus X_0$, then the generating function of the initial data can be written down as the sum
\[
\Phi(z) = \sum\limits_J \Phi_J (z) = \sum\limits_J \sum\limits_{\tau\in \Gamma_J} \Phi_{\tau,J}(z).
\]

\begin{theorem}\label{th_o_proizv_f}
The generating function $F(z)$ of the solution of the problem (\ref{razn_ur}),(\ref{nach_dan}) under the assumption (\ref{neravenstva}) and the generating function $\Phi(z)$ of the initial data are connected by the formula
\begin{align}\label{n-dim}
P(z) F(z) = \sum_J \sum_{\tau\in \Gamma_J} \Phi_{\tau, J}(z) P_\tau (z),
\end{align}
where $P_\tau(z)=\sum\limits_{\substack{\alpha\leqslant m \\ \alpha \nleqslant \tau}} c_\alpha z^\alpha$ are polynomials.
\end{theorem}

For $n=1$ it is easy to verify the statement of the Theorem \ref{th_o_proizv_f}. Indeed, in this case the generating function
\[
F(z)=\sum_{x=0}^{\infty} \frac{f(x)}{z^x}
\]
of the solution $f(x)$ to the initial value problem (\ref{n1_1}) with the coefficients \linebreak $(c_0, c_1, \ldots, c_m)$ and given initial data $(\varphi(0), \varphi(1), \ldots, \varphi(m-1))$ is a rational function:
\begin{equation}\label{1dimenstional}
F(z) = \frac {\sum\limits_{\alpha=1}^{m} \sum\limits_{x=0}^{\alpha-1} \frac{c_\alpha \varphi(x)}{z^{x-\alpha}} } { \sum\limits_{\alpha=0}^{m} c_\alpha z^\alpha}.
\end{equation}

For $n=2$ Theorem \ref{th_o_proizv_f} was proved in \cite{Lyapin2009} in connection with studying the rational Riordan arrays. For $n>1$ the proof was given in \cite{LeinLyapin2009}. The properties of generating function of solutions of a difference equation in rational cones of integer lattice were studied by T.~Nekrasova (see, e.g., \cite{Nekrasova}).

Theorem \ref{th_o_proizv_f} yields the following multidimensional analog of the De Moivre Theorem that is essential for the construction of algorithm:
\begin{theorem}\label{th_moivre}
The generating function $F(z)$ of the solution of the initial value problem (\ref{razn_ur}), (\ref{nach_dan}) under the assumption (\ref{neravenstva}) is rational if and only if the generating function $\Phi(z)$ of the initial data is rational.
\end{theorem}

The proof of Theorem \ref{th_moivre} was given in \cite{LeinLyapin2009}.

\section{Description of the Algorithm}

For $n=1$ the expression (\ref{1dimenstional}) consists of a finite number of terms, which makes the corresponding algorithm and computational procedure trivial. In this case the input data consists of two finite sets of numbers, namely: coefficients of the difference equation and the initial data. The output data of the algorithm is a rational function (\ref{1dimenstional}).

In the case when $n>1$, the initial data set $X_0$ is infinite. For $n=2$ the algorithm for computing the generating function $F(z)$ can be reduced to computation of a finite number of one-dimensional generating functions of sequences with elements along the coordinate axes. These elements are uniquely determined by coefficients of the corresponding one-dimensional difference equation and by finite set of the corresponding initial data (that is different for each sequence).

We now describe the structure of input data of the algorithm which consists of the matrices $c$, $C$ and $InData$.

The matrix $c = (c_{\alpha_1, \alpha_2})$, where $0\leqslant\alpha_1\leqslant m_1$, $0\leqslant\alpha_2\leqslant m_2$ has size $(m_1+1) \times (m_2+1)$. It's elements $c_{\alpha_1, \alpha_2}$ equal to the coefficients of two-dimensional difference equation (\ref{razn_ur}) if $(\alpha_1,
\alpha_2) \in A$ and equal 0 otherwise.

The coefficients of one-dimensional difference equations (that allow one to determine infinite initial data of two-dimensional initial value problem) are given by the matrix $C$. It consists of two columns: elements of the first column define starting point and direction of the corresponding one-dimensional initial data while elements of the second column are the coefficients of the one-dimensional difference equations, which define initial data on the corresponding direction.

We now consider the structure of the matrix $C$ more precisely. Suppose that for each $\xi_1 = 0, \ldots, m_1-1$ there exist numbers $\mu_{\xi_1} \in \mathbb N$ and  $c_{\mu_{\xi_1}}^{\xi_1}, c_{\mu_{\xi_1}-1}^{\xi_1}, \ldots, c_{0}^{\xi_1} \in \mathbb C$ such that a subset of the initial data $\{\varphi(\xi_1,x_2)\}_{x_2=0}^\infty$ along the horizontal axis satisfies one-dimensional initial value problem with the coefficients $(c_{0}^{\xi_1}, c_{1}^{\xi_1}, \ldots, c_{\mu_{\xi_1}}^{\xi_1})$ of the difference equation and with initial data $(\varphi (\xi_1,0), \varphi (\xi_1,1), \ldots, \varphi (\xi_1,\mu_{\xi_1}-1))$.
Then (\ref{1dimenstional}) implies that for each $\xi_1 = 0, \ldots, m_1-1$ the one-dimensional generating function of the sequence $\{\varphi(\xi_1,x_2)\}_{x_2=0}^\infty$ is rational.

Similarly, suppose that for each $\xi_2 = 0, \ldots, m_2-1$ there exist numbers $\nu_{\xi_2} \in \mathbb N$ and  $d_{\nu_{\xi_2}}^{\xi_2}, d_{\nu_{\xi_2}-1}^{\xi_2}, \ldots, d_{0}^{\xi_2} \in \mathbb C$ such that a subset of the initial data $\{\varphi(x_1,\xi_2)\}_{x_1=0}^\infty$ along the vertical axis are satisfies one-dimensional initial value problem with the coefficients $(d_{0}^{\xi_2}, d_{1}^{\xi_2}, \ldots, d_{\nu_{\xi_2}}^{\xi_2})$ of the difference equation and with initial data $(\varphi (0, \xi_2), \varphi (1,\xi_2), \ldots, \varphi (\nu_{\xi_1}-1,\xi_2))$.
Then (\ref{1dimenstional}) implies that for each $\xi_2 = 0, \ldots, m_2-1$ the one-dimensional generating function of the sequence $\{\varphi(\xi_1,x_2)\}_{x_2=0}^\infty$ is rational function. Both these functions can be easily computed.


Consequently, the matrix $C$ has the following structure
\begin{gather*}
C =
\begin{pmatrix}
  (1,0) & (c^0_{\mu_0}, c^0_{\mu_0-1}, \ldots, c^0_0) \\
  (2,0) & (c^1_{\mu_1}, c^1_{\mu_1-1}, \ldots, c^1_0) \\
  \ldots & \ldots \\
  (m_1,0) & (c^{m_1-1}_{\mu_{m_1-1}}, c^{m_1-1}_{\mu_{m_1-1}-1}, \ldots, c^{m_1-1}_0)\\
  (0,1) & (d^0_{\nu_0}, d^0_{\nu_0-1}, \ldots, d^0_0) \\
  (0,2) & (d^1_{\nu_1}, d^1_{\nu_1-1}, \ldots, d^1_0) \\
  \ldots & \ldots \\
  (0, m_2) & (d^{m_1-1}_{\nu_{m_1-1}}, d^{m_1-1}_{\nu_{m_1-1}-1}, \ldots, d^{m_1-1}_0)
\end{pmatrix}.
\end{gather*}

The initial data of one-dimensional difference equations are given by the matrix $InData$ of the size $M_1 \times M_2$, where $M_1 = \max \{\mu_{0}, \ldots, \mu_{m_1-1}\}$ and $M_2  = \max \{ \nu_{0}, \ldots, \nu_{m_2-1}\}$.

Note that some elements of $InData$ can be dependent from each other (or, in other words, some elements can be derived from another ones), so that it is necessary to define only a part of its' elements. We illustrate this fact on the following example.

Suppose that the subsets of the initial data of two-dimensional difference equation along the horizontal axis are given by the following three one-dimensional difference equations
\[
c_3^0\varphi(x+3,0)+c_2^0\varphi(x+2,0)+c_1^0\varphi(x+1,0)+c_0^0\varphi(x,0)=0,
\]
\[
c_2^1\varphi(x+2,1)+c_1^1\varphi(x+1,1)+c_0^1\varphi(x,1)=0,
\]
\[
c_2^2\varphi(x+2,2)+c_1^2\varphi(x+1,2)+c_0^2\varphi(x,2)=0.
\]
Similarly, let the subsets of the initial data along the vertical axis be given by the two one-dimensional difference equations
\[
d_2^0\varphi(0,y+2)+d_1^0\varphi(0,y+1)+d_0^0\varphi(0,y)=0,
\]
\[
d_3^1\varphi(1,y+3)+d_2^1\varphi(1,y+2)+d_1^1\varphi(1,y+1)+d_0^1\varphi(1,y)=0.
\]
Then the initial data of one-dimensional difference equations is given by a matrix
\[
InData =
\begin{pmatrix}
* & \varphi(1,2) & * \\ \varphi(0,1) & \varphi(1,1) & * \\ \varphi(0,0) & \varphi(1,0) & \varphi(2,0)
\end{pmatrix}.
\]
The computational procedure does not read elements denoted by $*$, but compute them from the other entries of the matrix and corresponding difference equations.

Consequently, the input data of the algorithm is finite.

Now we rewrite (\ref{n-dim}) in the form that is convenient for the implementation of the algorithm:
\begin{multline}\label{osn_form2}
F(z_1,z_2) =
 \left(\sum_{\xi_1=0}^{m_1-1} \sum_{\xi_2=0}^{m_2-1} \frac{P_{\xi_1,\xi_2}(z_1,z_2)}{z_1^{\xi_1+1}z_2^{\xi_2+1}} \varphi(\xi_1,\xi_2)\right. + \\
+ \sum_{\xi_1=0}^{m_1-1} \frac {P_{\xi_1,m_2}(z_1,z_2)}{z_1^{\xi_1+1}} \Phi_{m_2}^{\xi_1}(z_2) +   \\
+\left. \sum_{\xi_2=0}^{m_2-1} \frac {P_{m_1,\xi_2}(z_1,z_2)}{z_2^{\xi_2+1}} \Psi_{m_1}^{\xi_2}(z_1)\right) \slash {P(z)},
\end{multline}
where
\begin{equation} \label{osn_form2_1}
\Phi_{m_2}^{\xi_1}(z_2) = \frac {z^{m_2} \sum\limits_{\alpha=1}^{\mu_{\xi_1}} \sum\limits_{x_2=0}^{\alpha-1} \frac{c_\alpha^{\xi_1} \varphi(\xi_1,x_2+m_2)}{z_2^{x_2-\alpha}} } { \sum\limits_{\alpha=0}^{\mu_{\xi_1}} c_\alpha^{\xi_1} z_2^\alpha}
\end{equation}
and
\begin{equation} \label{osn_form2_2}
\Psi_{m_1}^{\xi_2}(z_1) = \frac {z^{m_1} \sum\limits_{\alpha=1}^{\nu_{\xi_2}} \sum\limits_{x_1=0}^{\alpha-1} \frac{d_\alpha^{\xi_2} \varphi(x_1+m_1,\xi_2)}{z_1^{x_1-\alpha}} } { \sum\limits_{\alpha=0}^{\nu_{\xi_2}} d_\alpha^{\xi_2} z_1^\alpha}.
\end{equation}

The right hand side of (\ref{osn_form2}) is split into 3 groups. The first group consists of a finite number of summands. The second and the third groups require computing of a finite number of one-dimensional generating functions (\ref{osn_form2_1}) and (\ref{osn_form2_2}).

The algorithm itself consists of three procedures that we present below.

The first procedure (Algorithm \ref{alg1}) computes the missing terms $\varphi(\xi_1,\xi_2)$ in the first sum of (\ref{osn_form2}). It is also used for computing generating functions in the second and third sums of (\ref{osn_form2}).

\begin{algorithm}
\caption{The algorithm for computing the initial data
$\varphi(x_0,y_0)$ at arbitrary point $(x_0,y_0)\in X_0$.}
\label{alg1}
\begin{algorithmic}
\Require Matrixes $c$, $C$, $InData$ and a point $(x_0,y_0)\in X_0$.
\Ensure A value of function of initial data $\varphi(x_0,y_0)$.
\Procedure{$\varphi$}{$c$, $C$, $InData$, $(x_0,y_0)$}

\State $m_1:=$ is a number of rows of the matrix $c$
\State $m_2:=$ is a number of columns of the matrix $c$

 \If{$x_0 < m_1$}
   \For{$i$ from $1$ to $m_1+m_2$}
      \If{$c[i,1] = (x_0+1,0)$}
        \State $c_0 := C[i,2]$ 
      \EndIf
   \EndFor
    \State $\varphi_0 := x_0$-th column of $InData$
    \State $l_0 := y_0$
 \EndIf

 \If{$y_0 < m_2$}
   \For{$i$ from $1$ to $m_1+m_2$}
      \If{$c[i,1] = (0,y_0+1)$}
        \State $c_0 := C[i,2]$ 
      \EndIf
   \EndFor
    \State $\varphi_0 := y_0$-th row of $InData$
    \State $l_0 := x_0$
 \EndIf
 \State $l := $ is a length of  $c_0$

    \For{$i$ from $l$ to $l_0$}
        \State $\varphi_0(i) := - \sum\limits_{j=0}^{l-1} \dfrac {c_0[j]}{c_0[m]}\varphi_0(i+j-l)$
    \EndFor

\Return $ \varphi_0(l_0)$
\EndProcedure
\end{algorithmic}
\end{algorithm}

\begin{algorithm}
\caption{The algorithm for computing one-dimensional generating function according to formulas (\ref{osn_form2_1}) or
(\ref{osn_form2_2}).} \label{alg2}

\begin{algorithmic}

\Require A vector $M$ of coefficients of a difference equation (with increasing order of indices); a vector $\Phi$ of initial data (with increasing order of indices) of the length equal to the length of $M - 1$; an integer $start$ (the initial value of summation); an integer $t\in\{1,2\}$; a variable $z_t$.
\Ensure Generating function, given by formulae (\ref{osn_form2_1}) or (\ref{osn_form2_2}) according to $t$.
\Procedure{$GF$}{$M$, $\Phi$, $start$, $t$, $z_t$}

\State $F:=0$
\State $n :=$  length of $M - 1$

\For {$\alpha$ from $1$ to $n+1$}
    \State $c_{\alpha-1} := M[\alpha]$
\EndFor

\For {$\alpha$ from $1$ to $n$}
    \State $\varphi[\alpha-1] := \Phi[\alpha]$
\EndFor

\For {$\alpha$ from $n$ to $start+n-1$}
    \State $\varphi[\alpha]:=0$
    \For {$j$ from $1$ to $n$}
        \State $\varphi[\alpha] := \varphi[\alpha] - \frac{c_{j-1}}{c_n} \cdot \varphi[\alpha- n+j-1]$
    \EndFor
\EndFor

\For {$\alpha$ from $0$ to $n-1$}
    \State $\varphi[\alpha] := \varphi[start+\alpha]$
\EndFor

\For {$\alpha$ from $1$ to $n$ }
    \For {$x$ from $0$ to $\alpha-1$}
        \State $F := F + \frac {c_\alpha \cdot \varphi[x]}{z_t^{x-\alpha+1}}$
    \EndFor
\EndFor

\State  $Q:=0$
\For {$\alpha$ from $0$ to $n$ }
    \State $Q := Q + c_\alpha \cdot z_t^\alpha $
\EndFor

\State $F:= F \cdot z_t^{start} / Q$

\Return $F$
\EndProcedure
\end{algorithmic}
\end{algorithm}

Algorithm \ref{alg3} is the main procedure that assemblies the two-dimensional generating function from one-dimensional generating functions obtained by Algorithm~\ref{alg2}.

\begin{algorithm}
\caption{The algorithm for computing two-dimensional generating function of the solution of the initial value problem (\ref{razn_ur}), (\ref{nach_dan}).}
\label{alg3}
\begin{algorithmic}[1]
\Require Matrices $c$, $C$, $InData$,
\Ensure Generating function $F(z)$ of the solution to the Cauchy problem.

\Procedure{$GenFunc$}{$c$, $C$, $InData$}

\State $m_1:=$ is a number of rows in/of the matrix $c$
\State $m_2:=$ is a number of cloumns in/of the matrix $c$

\State  $F:=0$

\For{$x_1$ from $1$ to $m_1-1$}
\For{$x_2$ from $1$ to $m_2-1$}
\State  $F:=F+\dfrac{P_{x_1,x_2}(z_1,z_2)}{z_1^{x_1} z_2^{x_2}} \cdot \varphi(c, C, InData, (x_0, y_0))$  
\EndFor
\EndFor

\For {$\xi_1$ from $0$ to $m_1-1$}
\State $\Phi := (\xi_1+1)$-th column of $InData$
\State $M := C[i,2]$ such that $C[i,1] = [\xi_1+1,0]$

\For {$i$ from $1$ to $\xi_1+\xi_2$}
    \If {$C[i,1]=(\xi_1+1,0)$}
        \State $M := C[i,2]$
    \EndIf
\EndFor

\State $F := F + \frac {P_{\xi_1,m_2}(z_1,z_2)}{z_1^{\xi_1+1}} \cdot GF(M,\Phi,m_2,2)$
\EndFor

\For {$\xi_2$ from $0$ to $m_2-1$}
\State $\Phi := (\xi_2+1)$-th row of $InData$
\For {$i$ from $1$ to $\xi_1+\xi_2$}
    \If {$C[i,1]=(0,\xi_2+1)$}
        \State $M := C[i,2]$
    \EndIf
\EndFor

\State $F := F + \frac {P_{m_1,\xi_2}(z_1,z_2)}{z_2^{\xi_2+1}} \cdot GF(M,\Phi,m_1,1)$
\EndFor

\State $F:=F/P(z_1, z_2)$

\State \textbf{return} $F$

\EndProcedure

\end{algorithmic}

\end{algorithm}

\section{Experiments}

The algorithm was implemented in Maple 2015. 
It was tested on the following examples.

\begin{example}
The binomial coefficients
\[
C_n^k=\frac{n!}{k!(n-k)!}
\]
give a solution to the Cauchy problem
\begin{align*}
f(x+1,y+1)-f(x,y+1)-f(x,y)=0
\end{align*}
with the initial data
\begin{align*}
\varphi(x,y)=\left\{
               \begin{array}{ll}
                 1, & \hbox{если $x\geqslant 0, y=0$;} \\
                 0, & \hbox{если $x=0, y \geqslant 0$.}
               \end{array}
             \right.
\end{align*}

In the example $m_1=m_2=1$ and the set of initial data
\[
X_0  = \{(x_1,x_2)\in \mathbb Z^2_+ \colon (x_1,x_2) \ngeqslant (1,1)\}.
\]
\begin{figure}[ht]
  \noindent \centering
  \includegraphics[width=70mm]{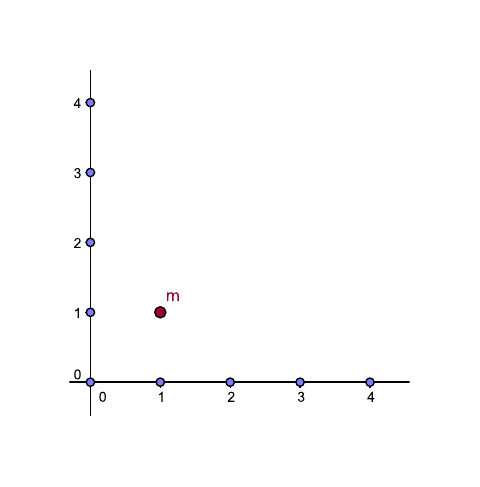}\\
  \caption{The initial data set $X_0$}
\end{figure}

The input data for the computational procedure:
\begin{gather*}
c =
\begin{pmatrix}
  -1 & 1 \\
  -1 & 0
\end{pmatrix},\\
InData =
\begin{pmatrix}
  1 & 1 \\
  0 & 1
\end{pmatrix},\\
C =
\begin{pmatrix}
  (0,1) & (-1,1) \\
  (1,0) & (1,0)
\end{pmatrix}.
\end{gather*}

The result is the generating function of the binomial coefficients:
\[
F (z,w) = \frac {1}{zw-w-1}.
\]
\end{example}

\begin{example}

Bloom studies the number of singles in all $2^x$ $x$-length bit
strings~\cite{bloom}, where a single is any isolated 1 or 0, i.e.
any run of length 1. Let $r(x,y)$ be the number of $x$-length bit
strings beginning with 0 and having $y$ singles. Apparently
$r(x,y)=0$ if $x<y$. Then $r(x,y)$ is:
\[
\begin{pmatrix}
        \vdots & \vdots & \vdots & \vdots & \vdots & \vdots & \vdots & \\
    0 & 0 & 0  & 0  & 0 & 0 & 1 & \ldots\\
    0 & 0 & 0  & 0  & 0 & 1 & 0 & \ldots\\
    0 & 0 & 0  & 0  & 1 & 0 & 5 & \ldots\\
    0 & 0 & 0  & 1  & 0 & 4 & 4 & \ldots\\
    0 & 0 & 1  & 0  & 3 & 3 & 9 & \ldots\\
    0 & 1 & 0  & 2  & 2 & 5 & 8 & \ldots\\
    1 & 0 & 1  & 1  & 2 & 3 & 5 & \ldots\\
\end{pmatrix},
\]
where the element $r(0,0)$ is in left lower corner.

In \cite{bloom} D.~Bloom proves that $r(x,y)$ is a solution to the Cauchy problem
\[
r(x+2, y+1) - r(x+1, y+1) - r(x+1, y) - r(x, y+1) + r(x,y) = 0
\]
with the initial data
\begin{align*}
&\varphi(0,0) = 1,\quad \varphi(1,0) = 0, \\
&\varphi(x,0) = \varphi(x-1,0)+\varphi(x-2,0),\ x\geqslant 2, \\
&\varphi(0,y) = 0,\ y\geqslant 1, \\
&\varphi(1,1)=1, \quad \varphi(1,y) = 0,\ y\geqslant 2.
\end{align*}

In this case the input data of the algorithm will be:
\begin{gather*}
c =
\begin{pmatrix}
  -1 & 1 \\
  -1 & -1\\
   0 & 1
\end{pmatrix},\\
InData =
\begin{pmatrix}
  1 & 0 \\
  0 & 1
\end{pmatrix},\\
C =
\begin{pmatrix}
  (0,1) & (-1,-1,1) \\
  (1,0) & (0,1)\\
  (2,0) & (0,0,1)
\end{pmatrix}.
\end{gather*}

The result is the generating function of the considered initial value problem:
\[
F (z,w) = \frac {z - 1}{z^2w - zw - w - z + 1}.
\]
\end{example}

\noindent \textbf{Acknowledgements.} The research was partially
supported by the state order of the Ministry of Education and
Science of the Russian Federation for Siberian Federal University
(task 1.1462.2014/K, first author) and by the Russian Foundation
for Basic Research, projects 15-31-20008-mol\_a\_ved,
15-01-00277-a (the first and the third author) and 14-01-00283-a
(first and second authors).

\end{document}